\documentclass[12pt]{amsart}
\usepackage{amscd}
\usepackage{verbatim}
\usepackage{amssymb, amsmath, amsthm, amscd,ifthen}
\usepackage[dvips]{graphics}
\usepackage[cp866]{inputenc}
\usepackage{graphicx}
\usepackage{epsfig}
\usepackage{xcolor}

\usepackage{amsmath,amssymb,amscd,}
\usepackage[all]{xy}

\usepackage{amssymb}

\newtheorem{thm}{Theorem}[section]

\newtheorem{Ex}[thm]{Example}

\newtheorem{cor}[thm]{Corollary}

\theoremstyle{definition}

\newtheorem{dfn}[thm]{Definition}

\usepackage{tikz}
\usepackage{tkz-euclide}
\usetikzlibrary{shapes,arrows,spy,positioning,snakes}

\title[Envy-free division]{Envy-free division \\ via configuration spaces}

\author[G. Panina]{Gaiane Panina}
\author[R. \v{Z}ivaljevi\'{c}]{Rade \v{Z}ivaljevi\'{c}}

\address[G. Panina]
{ St.\ Petersburg State University; St.\ Petersburg Department of Steklov Mathematical Institute}
\address[R. \v{Z}ivaljevi\'{c}]{Mathematical Institute of the Serbian Academy of Sciences and Arts (SASA), Belgrade}

\address{
  }

 \keywords{envy-free division, configuration space/test map scheme}

\begin{document}

\begin{abstract}
The classical approach to envy-free division and equilibrium problems arising in mathematical economics typically relies on Knaster-Kuratowski-Mazurkiewicz theorem, Sperner's lemma or some extension involving mapping degree. We propose a different and relatively novel approach where the emphasis is on configuration spaces and equivariant topology, originally developed for applications in discrete and computational geometry (Tverberg type problems, necklace splitting problem in the sense of  N.Alon and D. West, etc.). We illustrate the method by proving several relatives (extensions) of the classical envy-free division theorem of David Gale, where the emphasis is on preferences allowing the players to choose degenerate pieces of the cake.
\end{abstract}

 \maketitle \setcounter{section}{0}

\section{Overview and an informal introduction}\label{sec:Overview}

 Given some resource and a set of agents (players), one of the goals of
\emph{welfare economics} is to divide the resource among the agents in an \emph{envy-free} manner.
Envy-freeness, as a criterion of fair division, is the principle where every player feels that their share is at least as good as the share of any other agent, and thus no player feels envy.

\medskip
In a mathematical simplification a group of $r$  players want to divide among themselves a commodity, commonly referred to as the ``cake''. The simplest model of the cake is the interval $I=[0, 1]$, which should  be cut into $r$ pieces  by $r-1$ cuts.
Therefore a partition (cut) of the cake (in this model) is a sequence $x = (x_1,\dots, x_{r-1})$ where
\begin{equation}\label{eqn:cut}
    0 \leqslant x_1 \leqslant x_2 \leqslant \dots \leqslant x_{r-1} \leqslant 1 \, .
\end{equation}
The pieces of the cake arising from this division are the closed intervals (tiles)  $I_i = [x_{i-1}, x_i]\, (i=1,\dots, r)$, where $x_0 = 0$ and $x_r=1$.

\medskip
It would seem plausible that in actual applications all inequalities in (\ref{eqn:cut}) should be strict. Indeed the non-degeneracy of tiles is a salient feature of some classical results such as Gale's theorem (Theorem \ref{GaleThm}). However, as noticed by several authors in more recent publications \cite{AK, AK1, MeZe, S-H},  this condition may be too restrictive in some situations.

For example if the cake is partially burnt, the players are not hungry or there is some other reason which could make
some large pieces less desirable, then a degenerate piece may be preferred by some players.

One of our goals is to demonstrate that methods based on equivariant topology and new configuration spaces is well suited for the study of envy-free division  where \emph{degenerate tiles} may be preferred.

\subsection{Preferences} After a partition (\ref{eqn:cut}) is chosen, each of the players expresses his/her own individual \emph{preferences}, over the tiles arising from this particular cut, by pointing to one (or more) intervals which they like more than the rest.

The nature of the preferences may be unknown or hidden (black box preferences). For instance we do not assume the monotonicity of preferences, where a player satisfied by an interval $I$ is automatically satisfied with any $J\supseteq I$.

\medskip
Actually the preferences may be quite arbitrary, even irrational (from our point of view) as illustrated by the following example.

\begin{Ex} {\rm
Assume that two cuts are given
 $$C_1 = (1/4,1/3,2/3,3/4) \quad \mbox{\rm and}   \quad  C_2= (1/4,1/3,3/4,5/6)\, ,$$
where the two initial tiles $I_1 = [0,1/4]$ and $I_2= [1/4,1/3]$ are the same.
A player  may prefer the tile $I_1$ in the cut $C_1$, and the tile $I_2$ for the cut $C_2$.
So, although the tiles $I_1$ and $I_2$ are identical in both cuts, for some mysterious reason a player may prefer one over the other depending on the rest of the tiles.   }
\end{Ex}

\subsection{Conditions satisfied by all preferences} We emphasized that preferences can be quite arbitrary. There are however some conditions which are natural (intrinsic) or even unavoidable in the sense that they must be satisfied by all preferences.

\begin{Ex}{\rm \textbf{(Baguette's gorbushka example)} \label{Ex:gorbushka}
Suppose that each of the players prefers either the first or the last tile. This is a real-life preference in some countries (for instance Russia), where some people prefer the ``gorbushka'' (the ``heel'' or the end slice of the baguette).
If $r>2$  then envy-free division is clearly impossible since a loaf of bread has only two end pieces. }
\end{Ex}

The following two conditions may not be ``unavoidable''. However they are quite natural (reasonable) and they are not very restrictive so they are assumed throughout the paper as default conditions satisfied by all preferences.

\begin{itemize}
\item[$\mathbf{(P_{cl})}$] (The preferences are \emph{closed}) If a sequence of cuts converges to a limit cut and if in this sequence a player always prefers the tile labeled by $i$, then in the limit cut the player also prefers the tile labeled by  $i$.
\item[$\mathbf{(P_{cov})}$] (The preferences form a covering) Each player in each cut is expected to find an acceptable tile. In other words each player prefers at least one of the offered pieces.
\end{itemize}

\subsection{The role of degenerate tiles}

A classical theorem, found independently by Stromquist \cite{Strom} and by Woodall \cite{Wood} in 1980  states that, under the  conditions \emph{``no player prefers a tile which degenerates to a single point''}, it is possible to divide the cake
into $r$ connected tiles  and distribute these shares to  $r$ players in an\emph{ envy-free manner} (see Theorem \ref{GaleThm} for a formal statement). This result is referred to as Gale's theorem by some authors (after D. Gale who introduced important new ideas).

\medskip\noindent
\textbf{Remark}. The  conditions \emph{``no player prefers a tile which degenerates to a single point''} together with closeness condition implies, in particular, that ``{no player prefers an arbitrarily small tile}''.  That is, for a fixed  system of preferences, there exists $\varepsilon >0$ such that a tile shorter than $\varepsilon$ is never preferred.

\medskip
The ``gorbushka example'' (Example \ref{Ex:gorbushka}) shows that if the players are allowed to choose degenerate pieces then an envy-free division may not be possible without some restrictions.

\subsection{Preferences where the choice of degenerate tiles is permitted}
\label{sec:dte-pe}

In the present paper we want to understand how one can generalize or alter Gale's  theorem by allowing players to choose degenerate tiles. The following two types of conditions have already appeared in the literature.

\begin{itemize}
\item[$\mathbf{(P_{dte})}$] (All degenerate tiles are equal) In this type of preferences
all degenerate tiles are equal from the viewpoint of the players. Moreover each of them is informally equal  to the empty set.

\medskip\noindent
Example: Assume that   the cut-points are $(1/4,1/3,1/3,2/3,2/3)$. The condition $\mathbf{P_{dte}}$ says that  if a player prefers the degenerate segment $[1/3,1/3]$,  then  he also prefers the degenerate segment $[2/3,2/3]$.

\medskip\noindent
$\mathbf{P_{dte}}$ is the most important property used in the present paper. We shall see that this condition is essential for the construction of the new-style configuration spaces.

\item[$\mathbf{(P_{pe})}$] (Preferences satisfying partition equivalence) If two cuts produce one and the same collection of non-degenerate tiles, the preferences for these two cuts should define each other. That is,
 \begin{enumerate}
 \item A player prefers a non-degenerate tile for the first cut iff he  prefers the same tile for the second cut.
 \item A player prefers a degenerate tile for the first cut iff he prefers any of the degenerate tiles  for the second cut, and vice versa.
 \end{enumerate}

 Note that $\mathbf{P_{dte}}$ is a consequence  of (2).

\medskip\noindent
Example: Assume that the cut-points are $(1/4,1/3,1/3,2/3)$. The  non-degenerate tiles are $[0,1/4], \  [1/4,1/3], \ [1/3,2/3], [2/3,1]$. The following cuts  $(1/4,1/3,2/3,2/3)$ and $(1/4,1/3,2/3,1)$  are formally different although they have the    same non-degenerate tiles.
The $\mathbf{(P_{pe})}$ condition states that if a player prefers the segment $[1/4,1/3]$ in the first cut, he keeps preferring the segment $[1/4,1/3]$ for the other cut.
\medskip

The condition $\mathbf{P_{pe}}$  is used only once in the paper (in Theorem \ref{thm:AK4.1}).

\medskip
An easy exercise shows that with the Baguette's gorbushka preferences the $\mathbf{P_{dte}}$ condition makes fair division possible. Indeed, if the ``gorbushka'' degenerates to a point, it becomes equivalent to any other degenerate segment. So if the cake master makes all the cuts at the end-point $0$, and so creates $r-1$ degenerate tiles, every player becomes happy.
\end{itemize}

\subsection{What kind of divisions are allowed}\label{sec:division}

  A division of the cake is called \emph{full} if each non-degenerate tile is allocated to one of the players, or in other words \emph{no tiles are dropped}.

Other divisions may be possible. For example some non-degenerate tiles may be dropped because some players prefer not to choose any of the available non-degenerate tiles (we say that they prefer an empty piece).

Sometimes it may be  necessary to allow exceptions (see Theorem \ref{ThmPrime}) when one player (in a very exceptional situation) is allowed to have two pieces of the cake.

\subsection{\textbf{What is  known }}
\begin{enumerate}
  \item If $r$ is a prime power then a full, envy-free division always exists if the preferences satisfy the $\mathbf{P_{pe}}$ condition. See \cite{AK} and \cite{AK1}, and Theorem \ref{ThmKA} in the present  paper.
This result of Avvakumov and Karasev \cite[Theorem 3.3]{AK} was established earlier by Meunier and Zerbib \cite[Theorem 1]{MeZe} in the cases when $r$ is a prime or $r=4$, see also Segal-Halevi \cite{S-H} where the result was conjectured and proved in the case $r=3$.
  \item As shown in \cite{AK} if $r$ is not a prime power the theorem of Avvakumov and Karasev is no longer true.
  \item Assume that we are allowed to  drop some of the tiles. Then without restrictions on  $r$  a not necessarily full envy-free  division always exists under the condition $\mathbf{P_{dte}}$, see Theorem \ref{DropThm} and  \cite{AK, MeZe} for a proof.
\end{enumerate}

\subsection{\textbf{What is new in the paper}}
\begin{enumerate}
\item We modify the rules of the game. As before there are $r-1$ cuts which partition the interval $[0,1]$ into $r$ tiles   (possibly degenerate). In the presence of degenerate tiles only one of the guests may be given two tiles (instead of one) in which case one of them must be the rightmost tile (the right gorbushka).
Then, if either $r$ is a prime number or $r=4$, under the condition $\mathbf{P_{dte}}$, a full envy-free division always exists with no other restrictions (Theorem \ref{ThmPrime}).
\item If $r$ is neither a prime nor equal to $4$, the statement above  is no longer true (Theorem \ref{ThmNonPrime}).
 \item  Let $r$ be a prime number or $r=4$.   Assume that in the division we are allowed to drop \emph{at most one  of the non-degenerate tiles} and the condition  ``{no player prefers an arbitrarily small tile number $r$, unless it degenerates to a point}''.  That is, for a fixed  system of preferences, there exists $\varepsilon >0$ such that a non-degenerate tile $r$ shorter than $\varepsilon$ is never preferred. Then  an envy-free  division always exist  under the condition $\mathbf{P_{dte}}$ (Corollary \ref{CorDrop}).
 \item We give an alternative proof of the result of Avvakumov and Karasev: If $r$ is a prime power then an envy-free full partition always exist under the condition  $\mathbf{P_{pe}}$.
\item Consider  $r-1$ cuts which partition the interval $[0,1]$ into $r$ possibly degenerate tiles. We assume that the tiles may be distributed among the players without any restrictions, in particular a player may receive several tiles at a time.
Then under the condition $\mathbf{P_{dte}}$, a full envy-free division always exists if  $r=p^a$ is a prime power (Theorem \ref{ThmPrimePower2}).
\end{enumerate}

 \subsection{\textbf{New narrative and the novelty of the method.}}
 We approach the cake division problem from a slightly different point of view by offering a new  narrative which is, in our opinion, helpful for the  formalization and construction of configuration spaces.

 We assume that the cake division supervisor (cake master) cuts the cake into tiles, and then puts them on display in  (transparent) cake boxes, prepared in advance and numbered from $1$ to $r$. This is called a partition/allocation of the cake because it involves both a cut (partition) of the cake and (an arbitrary) allocation of the tiles in boxes.

 After that the players approach the boxes and point to some of them, expressing their preferences.

 The condition $\mathbf{P_{dte}}$ allows the cake master to ignore  the degenerate tiles, meaning that they are not allocated at all. Instead, the degenerate tiles are put in the trash, and as a result after a partition/allocation some of the boxes may  be empty.

 A quite natural assumption in this setting is  the \emph{equivariance of the preferences}. It says that the players pay attention only on the content, disregarding the numbers (labels) on the boxes. In other words each player prefers one and the same piece of the cake,  ignoring the numbering of the box where it is placed in the current allocation.

 An example of a non-equivariant preference arises if a (superstitious) thief always avoids the  box no.\ 13,
regardless of its  content.

 \bigskip

\section{Comparison of the old and the new setting}

The classical or the ``old'' approach to envy-free division typically relies on Knaster-Kuratowski-Mazurkiewicz theorem, Sperner's lemma, Brouwer fixed-point theorem, or some extension involving mapping degree. Our approach, where the emphasis is on configuration spaces and equivariant topology, is (informally) referred to as the ``new'' approach.

\subsection{Configuration spaces, old and new}
\label{sec:conf-old-new}

A \emph{cut} of the segment $I=[0,1]$ by $r-1$ cut-points is (as in Section \ref{sec:Overview}) recorded  as  $x=(x_1,...,x_{r-1})\in I^{r-1}$, where the individual cut-points $x_i$ are listed in the increasing order. The tiles $I_i = [x_{i-1}, x_i]$ are automatically labeled by $[r]=\{1,2,...,r\}$, also from left to the right.
The cut $x$ can be uniquely reconstructed from the vector $z = (z_1,\dots, z_r) \in \Delta^{r-1} := {\rm conv}\{e_i\}_{i=1}^r\subset \mathbb{R}^r$,  where $z_i\in [0,1]$ is the length of the tile $I_i$.

Therefore, in the traditional setting, the \emph{configuration space of all possible cuts} is the standard simplex $\Delta^{r-1}$. Its boundary corresponds to degenerations of some of the tiles.

\medskip

Assuming that there are $r$ different boxes, also labeled by $[r]=\{1,2,...,r\}$,
 let us consider an \emph{allocation map} $\alpha : [r] \rightarrow [r]$  whose function is to distribute the tiles in   the boxes. To this end, there are no restrictions on $\alpha$: some of the boxes may contain more than one tile.

\begin{dfn}
Two allocations with one and the same cut are called \emph{equivalent} if they differ only by allocations of degenerate tiles.
\end{dfn}
This definition formalizes (on the level of configuration spaces) the property $\mathbf{P_{dte}}$ which says that ``all degenerate tiles are equal from the viewpoint of the players''.

\medskip\noindent
{Example:} Assume that  the cut-points are $(1/4,1/3,1/3,2/3,2/3)$, which means that the tiles $I_3$ and $I_5$ are degenerate.  All tiles $I_i$, labeled by $i\in [6]$, should be  allocated to six boxes $B_j$, labeled by $j\in [6]$.
The allocation
\newline $B_1=\{1,2,3\},B_2=\emptyset, B_3=\emptyset, B_4=\{4\},B_5=\{5\}, B_6=\{6\}$ is equivalent to
\newline $B_1=\{1,2\},B_2=\emptyset, B_3=\emptyset, B_4=\{4\},B_5=\{5,3\}, B_6=\{6\}$, and also to
\newline $B_1=\{1,2\},B_2=\{3\}, B_3=\emptyset, B_4=\{4\},B_5=\{5\}, B_6=\{6\}$.

\medskip

It is instructive to think that for an equivalence class only non-degenerate tiles are allocated, whereas degenerate ones are put in the trash.
So the equivalence class of the above three  allocation is recorded as \
\newline $B_1=\{1,2\},B_2=\emptyset, B_3=\emptyset, B_4=\{4\},B_5=\emptyset ,B_6=\{6\}$.

\medskip

 The space of all possible equivalence classes of  cuts and allocations forms a configuration space. This configuration space can be described as a simplicial complex obtained by gluing together different copies $\Delta_\alpha$ of the standard simplex $\Delta^{r-1}$, one copy for each allocation function $\alpha : [6] \rightarrow [6]$.

 \medskip
This construction has many variations and modifications, for example one can put some restrictions on the mapping $\alpha$. Different restrictions yield different configuration spaces and the design of a proper configuration space is guided and motivated by the envy-free division problem we are interested in.

All these configuration spaces form a \emph{species} which is closely tied to the class of \emph{chessboard complexes} (Section \ref{sec:chessboard}).

 \medskip

 Summarizing the main difference between the new and the old approach we observe that the configuration spaces are respectively of the following two types:

 \begin{itemize}
   \item The old-style configuration space is the simplex $\Delta^{r-1}$, whose points record cuts of the interval $I=[0,1]$.
   \item Characteristic for the new approach is a configuration space $\mathcal{C}$, whose points record both the cuts of $I$ and allocations of non-degenerate tiles into boxes.
 \end{itemize}

 \subsection{Preferences, old and new}
 \label{sec:prefs-old-new}

\begin{dfn}
The old-style preferences of $r$ players is a matrix of subsets $(A^j_i)_{i,j=1}^r$ of the standard simplex $\Delta^{r-1}$. The subsets are interpreted as preferences in the usual sense as follows:
\begin{equation}\label{eqn:prefs-1}
x\in A^j_i \ \ \Leftrightarrow  \, \mbox{ {\rm in the cut} } x \mbox{ {\rm the player} } j \mbox{ {\rm prefers the tile} } i \, .
\end{equation}
The new-style preferences of $r$ players is a matrix of subsets $(B^j_i)_{i,j=1}^r$ of the configuration space $\mathcal{C}$. Analogously, they are interpreted as:
\begin{equation}\label{eqn:prefs-2}
\begin{split}
(x, \alpha) \in B^j_i \ \ \Leftrightarrow  \, &  \mbox{ {\rm in the cut} } x \mbox{ {\rm and the allocation} }
\alpha \mbox{ {\rm the player} } j\\
& \mbox{ {\rm prefers the content of the box } } i \, .
\end{split}
\end{equation}
\end{dfn}

We assume throughout the paper that:
\begin{itemize}
\item[$\mathbf{(P_{cl})}$]  The preferences are \emph{closed}, meaning that they  are represented by closed subsets of the configuration space.
\item[$\mathbf{(P_{cov})}$] The preferences are \emph{covering}, in the sense that for each $j$ the union of the sets
$A_i^j,\, (i=1,\dots, r)$ is equal to the whole configuration space.
\end{itemize}

\medskip
The symmetric group $S_r$  acts   on $\mathcal{C}$ by renumbering the boxes:
    $$\sigma(x,\alpha):= (x,\sigma \circ \alpha)\, .$$

 We assume throughout the paper that the new-style preferences are \emph{equivariant} in the sense that for the action of each $\sigma\in S_r$,
\begin{equation}\label{eqn:equivariance-1}
  (x, \alpha ) \in B_{i}^j  \Leftrightarrow  \sigma(x, \alpha )\in B_{\sigma(i)}^j   \, .
\end{equation}
The justification of (\ref{eqn:equivariance-1}) is quite obvious. In light of (\ref{eqn:prefs-2}) $(x, \sigma\circ \alpha) \in B^j_{\sigma(i)}$ if and only if in the cut $x$ and the allocation $\sigma\circ \alpha$ the player $j$
prefers the content of the box $\sigma(i)$. Since $(\sigma\circ\alpha)^{-1}(\sigma(i)) = \alpha^{-1}(i)$ the condition (\ref{eqn:equivariance-1}) expresses the idea that players
make their decisions solely on the content of the boxes, not on their current labels.

\subsection{Preferences $\mathbf{P_{pe}}$}\label{sec:formal-dte=pe} Here we formalize the  property  $\mathbf{P_{pe}}$ ''The preferences are determined by non-degenerate tiles'' (introduced in Section \ref{sec:dte-pe}). For this purpose we provide a dictionary of some of the relevant terms and definitions (used here and elsewhere in the paper).

\bigskip\noindent
Degenerate tiles  for a cut $x\in \Delta^{r-1}$: \quad $Deg(x) = \{i\in [r] \mid \vert I_i\vert  = x_i- x_{i-1} = 0\}$

\medskip\noindent
Essential tiles for a cut $x$: \quad $Ess(x) = \{I_i \mid x_i - x_{i-1}> 0\}$

\medskip\noindent
Labeling of essential tiles: \quad For $I\in Ess(x) \quad \nu(I,x) = i \Leftrightarrow I = I_i = I_i(x)$

\medskip\noindent
Partition equivalent cuts:  \quad $x\sim_{pe} x' \Leftrightarrow Ess(x)=Ess(x') $

\medskip\noindent
We say that  the preferences $(A_i^j)$  are  $\mathbf{(P_{pe})}$ if whenever $x\sim_{pe} x'$, we have
\begin{enumerate}
    \item  $ \quad x\in A_{\nu(I,x)}^j
\Leftrightarrow x'\in A_{\nu(I,x')}^j
$, and
      \item  $
 (\exists i_1\in Deg(x)) \, x\in A_{i_1}^j   \, \, \Rightarrow \, \, (\forall i_2\in Deg(x')) \, x'\in A_{i_2}^j
.$
                                                                   \end{enumerate}

The {partition equivalence} condition is closely related to {\em pseudo-equivariance assumptions} from \cite[Section 4.1]{AK}.
Indeed, $x$ and $x'$ are in the same partition equivalence class if and only if one can be obtained from the other by an identification of the form $\sigma_{FGZ}$, described in \cite{AK}.

\subsection{Envy-free division, old and new}
\begin{dfn} Let $(A^j_i)_{i,j=1}^r$ be a matrix of old-style preferences.
A point $x\in \Delta^{r-1}$ describes an\emph{ envy-free division} for preferences
$(A^j_i)_{i,j=1}^r$ if there exists a permutation $\pi\in S_r$ such that
$$
    x \in \bigcap_{j=1}^r A_{\pi(j)}^j \, .
$$
\end{dfn}

Here is the counterpart of this definition for the new-style preferences.

\begin{dfn} Let $(B^j_i)_{i,j=1}^r$ be a matrix of new-style preferences.
A point $\hat{x} = (x,\alpha) \in \mathcal{C}$  provides an\emph{ envy-free division} for preferences  $(B^j_i)_{i,j=1}^r$ if
\begin{equation}\label{eqn:non-emty-envy-free}
(x,\alpha)\in \bigcap_{j=1}^r B_j^j  \, .
\end{equation}
\end{dfn}
By the equivariance of preferences (\ref{eqn:equivariance-1}) we immediately observe that if an envy-free division $(x,\alpha)$ exists then for each  permutation $\pi\in S_r$
\begin{equation}\label{eqn:pi-envy-free}
 (x,\pi\circ\alpha) \in \bigcap_{j=1}^r B_{\pi(j)}^j \, .
\end{equation}
For the same reason if (\ref{eqn:pi-envy-free}) holds for one permutation than it holds for all of them and  (\ref{eqn:non-emty-envy-free}) is an immediate consequence.

\bigskip
In the terminology of this section the envy-free division theorems of Gale and Avvakumov-Karasev can be formulated as follows.

\begin{thm}\label{GaleThm}  {\rm (Gale's theorem \cite{G})}
Let $(A_i^j)_{i,j=1}^r, \ A_i^j\subset \Delta^{r-1}$ be a matrix of the old-style closed covering preferences.
Assume that nobody prefers a degenerate tile, that is, $A_i^j\cap \Delta_i^{r-1} = \emptyset$ for each $i$ and $j$, where $\Delta_i^{r-1}$ is the face of $\Delta^{r-1}$ opposite to the vertex labeled by $i$.
Then there exists an envy-free division, that is, a permutation $\sigma\in S_r$ such that
\[
     \bigcap_{i\in [r]}  A_{\sigma(j)}^{j} \neq\emptyset \, .
\]
\end{thm}

\begin{thm}\label{ThmKA}{\rm (\cite[Theorem 4.1]{AK})}\label{thm:AK4.1}
Assume that  $r$  is a prime power. Assume we have the old-style closed, covering  preferences  $(A_i^j)_{i,j=1}^r, \ \ A_i^j \subset \Delta^{r-1}$ satisfying  the $\mathbf{P_{pe}}$ condition.

Then there exists a permutation $\pi\in S_r$ such that
\[
    \bigcap_{j\in [r]} A^j_{\pi(j)} \neq \emptyset \, .
\]
In other words there exists a partition of $[0,1]$ into at most $r$ non-degenerate intervals such that  each non-degenerate interval is given to a different player,
the remaining players are not given anything (they are given ``empty pieces''), and this distribution is envy-free from the viewpoint of each of the players.

\end{thm}

\begin{thm}\label{DropThm}{\rm   (\cite{AK})}
  Let $(A_i^j)_{i,j=1}^r$ be the old-style closed, covering preferences.
  Degenerate tiles may be preferred but the condition $\mathbf{P_{dte}}$ is satisfied.
 Then there exists an envy-free division where some tiles may be dropped.
\end{thm}

\section{Configuration space $\Delta_{r,r-1}\ast [r]$} \label{Sec1} In this section we consider  the configuration space of partitions and allocations with the following restrictions    on the allocation function $\alpha$:
 \begin{enumerate}
   \item  At most two non-degenerate tiles are allocated to the same box.
   \item If two non-degenerate tiles appear in the same box than one of the tiles is necessarily the last tile $r$ (the right ``gorbushka'' tile).
 \end{enumerate}

So, the tiles are allocated one tile per box with the unique exception: the tile $r$ can be put with some other tile in one and the same box.

\medskip
Denote the associated configuration space by $\mathcal{C}_1$. This space, which appeared in the literature in the framework of the \emph{colored Tverberg theorem},  is known as the \emph{chessboard complex}  $\Delta_{r,r-1}*[r]$. We need some special topological properties of this space, however the reader can safely take them for granted (at least for the first reading). For a more detailed exposition and the original references we refer the reader to the Appendix (Section \ref{sec:appendix}).

\medskip

The configuration space $\mathcal{C}_1$ is used in the proof of the following Gale-type theorem which we informally call the ``right gorbushka'' theorem.

\begin{thm}\label{ThmPrime}
Suppose that either $r$ is a prime number or $r=4$.
Let $$(B_i^j)_{i,j=1}^r,\ \  B_i^j\subseteq \mathcal{C}_1$$ be a matrix of the new-style preferences which are \emph{closed, covering and equivariant}. Then there exists an envy-free division:
\begin{equation}\label{eqn:envy-free-non-empty-2}
     \bigcap_{i\in [r]}  B_i^i \neq\emptyset \, .
\end{equation}
In other words, there exists a cut and an allocation (of the type described above) such that the player $j$ prefers the content of the box $j$ for all $j$. More explicitly, each player obtains not more than one non-degenerate tile, possibly with only one  exception, when a player gets two tiles. In that case one of these tiles is labeled by $r$.
\end{thm}
\proof
Assuming the contrary, suppose that the intersection (\ref{eqn:envy-free-non-empty-2})  is empty. It follows from the equivariance of the preferences that the intersection (\ref{eqn:non-emty-envy-free}) is also empty for each permutation $\pi\in S_r$.

\medskip
Let us construct an \emph{equivariant} test map $F : \mathcal{C}_1 \rightarrow \mathbb{R}^r$ which records the preferences of the players.
We begin by replacing closed preferences $B_i^j$ by slightly larger open equivariant preferences $O_i^j$, so that the corresponding intersection is also empty for each permutation $\sigma\in S_r$
\begin{equation}\label{eqn:open-empty}
     \bigcap_{j\in [r]}  O_{\sigma(j)}^j =\emptyset \, .
\end{equation}

Next, choose  (for each $j\in [r]$) an equivariant partition of unity $\{f_i^j\}_{i=1}^r$ subordinated to the cover $\{O_i^j\}_{i=1}^r$.
Since the preferences are equivariant, we can  make these  functions also equivariant:
\begin{equation}\label{eqn:equivariance}
f_{\sigma(i)}(\sigma(x, \alpha)) = f_{i}^j(x, \alpha) \, .
\end{equation}

By construction  $f_i^j(x, \alpha) > 0$ implies $(x, \alpha)\in O_i^j$.

\medskip
By averaging $F_i=\frac{1}{r}\sum_{j=1}^{r}f^j_i$ we obtain a vector-valued $S_r$-equivariant function  $$F=(F_1,F_2,\dots,F_r):\mathcal{C}_1 \rightarrow \mathbb{R}^r.$$
 Let $\widehat{F}  : \mathcal{C}_1 \rightarrow \mathbb{R}^r/D$ be the (equivariant) map obtained by composing the map $F$ with the projection $\mathbb{R}^r \rightarrow \mathbb{R}^r/D_{\mathbb{R}^r}$, where $D_{V^r}$ stands for the diagonal $D_{V^r} \subset V^r$. We claim that the map $\widehat{F}$ must have a zero.

\medskip
Otherwise there arises an equivariant map $\mathcal{C}_1 \rightarrow S^{r-2}$. However this is a contradiction, see  \cite[Proposition 4.2]{BMZ}  and \cite[Proposition 3]{vz11}.
Note that in these papers the configuration space $\mathcal{C}_1$ served as the key building block for the configuration space used in the colored Tverberg theorem.

\medskip

It remains to apply D. Gale's original argument \cite{G} (see also \cite{AK}). We have proven the existence of $(x, \alpha)\in \mathcal{C}_1$ such that the matrix $(f_i^j(x, \alpha))$ is doubly stochastic. By the Birkhoff-von Neumann theorem it lies in the  convex hull of
permutation matrices.

It follows that there exists a permutation $\sigma \in S_r$ such
that $f_\sigma(j)^{j}(x, \alpha)>0$ for all $j$. As a consequence the intersection (\ref{eqn:open-empty}) is also empty which is a contradiction.   \qed

\bigskip

Here is a corollary which refines Theorem \ref{DropThm} about not necessarily full divisions, which means that some of the tiles may not be distributed at all.

\begin{cor}\label{CorDrop} Assume that $r$ is prime or $r=4$, and let $(A_i^j)_{i,j=1}^r \, (A_i^j \subseteq \Delta^{r-1})$ be a matrix  of old-style (covering, closed, $\mathbf{P_{dte}}$) preferences.
We do not exclude that a degenerate tile is preferred by some players. We also assume that  there exists $\varepsilon >0$ such that a non-degenerate tile $r$ shorter than $\varepsilon$ is never preferred.

Then there always exists an envy-free division such that each player gets a preferred tile (or nothing, if he prefers so),
and at most one tile is dropped.
\end{cor}

\proof

Given $(A_i^j)$,  let us create new-style preferences $(B_i^j)=\psi (A_i^j), \  B_i^j\subseteq \mathcal{C}_1$ by the following rule.

Imagine you are a player $j$, the segment $I$ is cut into $r$ tiles by a cut $x$, and you know which tiles you prefer according to $(A_i^j)_{i=1}^r$. These tiles are allocated to $r$ boxes by an allocation function $\alpha : [r]\rightarrow [r]$, and you need to specify which boxes you prefer.

\bigskip

An empty box is preferred whenever a degenerate tile is preferred according to $(A_i^j)$. A non-empty box is preferred whenever one of the tiles
contained in the box is preferred according to $(A_i^j)$. We add one exceptional rule: If the box contains the tile $r$ and no other tiles, the box is preferred not only if the tile $r$ is preferred, but also if a degenerate tile is  preferred according to $(A_i^j)$.

\medskip
More formally,

\begin{enumerate}

  \item If the box $i$ is empty, then $(x, \alpha)\in B_i^j$   iff
   \begin{enumerate}
     \item There is a degenerate tile $m$ in the cut $x$, and
     \item $x\in A_m^j$.
   \end{enumerate}

  \item If the box $i$ has exactly one non-degenerate tile $m\neq r$,   then
  $$(x,\alpha)\in B^j_i \Leftrightarrow x\in A_m^j.$$

\item If the box $i$ contains two non-degenerate tiles, $r$ and some $m\neq r$,   then
  $$(x,\alpha)\in B^j_i \Leftrightarrow  \hbox{either \ \ } x\in A_r^j  \hbox{\ or\ } x\in A_m^j.$$

 \item If the box $i$ has exactly one non-degenerate tile $ r$,   then

 \begin{enumerate}
   \item If there are no degenerate tiles in the cut $x$,
    $$(x,\alpha)\in B^j_i \Leftrightarrow x\in A_r^j.$$

   \item If there is a degenerate tile $m$ in the cut $x$,
    $$(x,\alpha)\in B^j_i \Leftrightarrow  \hbox{either \ \ } x\in A_r^j \hbox{\ or\ } x\in A_m^j.$$

 \end{enumerate}
\end{enumerate}

\bigskip

The new-style preferences $(B_i^j)$ are clearly equivariant and covering. They are also closed -- this is where the condition ``short tile $r$ is never preferred unless it degenerates to a point'' is used. Let us briefly comment on how closeness fails without this condition. Assume there is a sequence of cuts $(x^{(n)})_{n\in \mathbb{N}}$  such that the length of the tile $r$ tends to zero, and the tile $r$ is preferred according to $(A_i^j)$.  Consider an allocation function $\alpha_k,\  (k=1,\dots,r-1)$ which places  the tile $r$ in the same box with the tile $k$, assuming that the tile $k$ is non-degenerate. For $(x^{(n)}, \alpha_k)$ the box with the tiles $r$ and $k$ is preferred. Therefore, in the limit, when $r$ degenerates, the box with the tile $k$ must be preferred, which is leads to a contradiction.

\medskip\noindent
\textbf{Remark.} For the above defined preferences $(B_i^j)= \psi (A_i^j)$, if there are no degenerate tiles, an empty box may exist, but an empty box is never preferred.

\bigskip

Now let us apply Theorem \ref{ThmPrime} to $(B_i^j)=\psi (A_i^j)$. There exists an envy-free division such that each player $j$ prefers the box $j$  for all $j$.
If the box $j$ is empty, the player $j$ prefers a degenerate tile. It exists by the above remark, so let the player $j$ take the degenerate tile.

If the box $j$ contains a unique tile which is not $r$, then the player prefers the tile.

If the box $j$ contains two tiles (this might happen for one box only), one of these tiles is preferred. Drop the other one.

If the box $j$ contain   a unique tile which is  $r$, then the player either prefers the tile (then there is no need to drop the tile), or he prefers the degenerate tile. In this case drop the tile $r$.\qed

\begin{thm}\label{ThmNonPrime}
Assume that $r$  is neither a prime number nor equal to $4$.
 Then there exist $(B_i^j)_{i,j=1}^r,  \, B_i^j \subseteq \mathcal{C}_1$,  the new-style (closed, covering, equivariant) preferences
 such that
\[
     \bigcap_{i\in [r]}  B_i^i =\emptyset \, .
\]

In other words, there exists no envy-free division if empty pieces may be preferred, and one of the players may be satisfied by getting a disconnected share consisting of two tiles one of which is $r$.
\end{thm}

\proof  Construct an equivariant function
$F=(F_1,F_2,\dots,F_r):\mathcal{C}_1 \rightarrow \mathbb{R}^r$  whose image belongs to the hyperplane $x_1+\dots +x_r=1$ and misses the diagonal $x_1=\dots =x_r$. If $r$ is not a prime the existence of such  function follows from \cite[Proposition 4.2]{BMZ}.

Let $\widetilde{B}_i^j$ be the preferences not depending on $j$ (meaning that all the players have the same preferences) where
$$\widetilde{B}_i^j=\widetilde{B}_i:=F^{-1}\{x_i\neq 0\}\, .$$

They are covering and equivariant, however they are open, rather than closed (as required). Still $\widetilde{B}_i$ can be replaced by a slightly smaller closed sets $B_i$ with all the necessary properties.
It remains to observe that $\bigcap_{i\in [r]}  B_i^i =\emptyset$. \qed

\section{Configuration space $\Delta_{r,2r-1}$}\label{Sec2}

In this section we use our method to give a new proof of the envy-free theorem of Avvakumov and Karasev \cite[Theorem 4.1]{AK}. Recall that this result was established earlier by Meunier and Zerbib \cite[Theorem 1]{MeZe} in the cases when $r$ is a prime or $r=4$, see also Segal-Halevi \cite{S-H} where the result was conjectured and proved in the case $r=3$.

\medskip

We modify the configuration space. Let us consider $2r-2$ cuts of the interval $I=[0,1]$, so altogether there are  $2r-1$ tiles. Allocate the tiles to $r$ boxes by \emph{``admissible''} allocation functions $\alpha : [2r-1]\rightarrow [r]$,  permitting at most one non-degenerate tile in each of the boxes.

\medskip
As a consequence, in the associated configuration space $\mathcal{C}_2$ only those cuts that create at least $r-1$ degenerate tiles are permitted.

The configuration space  $\mathcal{C}_2$ is again a chessboard complex $\Delta_{r,2r-1}$, see the Appendix (Section \ref{sec:appendix}). Its elements are (equivalent classes of) pairs $(x,\alpha)$ where $x$ is a cut with at least $r-1$ degenerate tiles and $\alpha$ is an associated admissible allocation function.

\begin{thm}\label{ThmPrimePower}
Let $r$ be a prime power.
  Let $$(C_i^j)_{i,j=1}^r,\ \  C_i^j\subseteq \mathcal{C}_2$$ be a matrix of the new-style (closed, covering, equivariant) preferences. Then there exists an envy-free division:
\begin{equation}\label{eqn:empty-4}
     \bigcap_{i\in [r]}  C_i^i \neq\emptyset \, .
\end{equation}
In other words there exists a cut and an \emph{admissible allocation} of non-degenerate tiles into the boxes, such that the player $j$ prefers the box $j$ for all $ j$.
\end{thm}
\proof

The proof follows the pattern of the proof of Theorem \ref{ThmPrime}.
Assume, for the sake of contradiction, that the intersection (\ref{eqn:empty-4}) is empty.
As before we construct an (equivariant) test map $F : \mathcal{C}_2 \rightarrow \mathbb{R}^r$ which records the information provided by the preferences $C_i^j$.

Initially we replace $C_i^j $ with slightly larger open preferences $O_i^j$ and choose  (for each $j\in [r]$) an equivariant partition of unity $\{f_i^j\}_{i=1}^r$ subordinated to the cover $\{O_i^j\}_{i=1}^r$.
\medskip
Since the preferences are equivariant, we can  make these  functions also equivariant:
\begin{equation}\label{eqn:equivariance-3}
f_{\sigma(i)}^j(\sigma(x, \alpha)) = f_{i}^j(x, \alpha) \, .
\end{equation}
By averaging $F_i=\frac{1}{r}\sum_{j=1}^{r}f^j_i$ we obtain a vector-valued $S_r$-equivariant function  $$F=(F_1,F_2,\dots,F_r):\mathcal{C}_1 \rightarrow \mathbb{R}^r.$$
 Let $\widehat{F}  : \mathcal{C}_1 \rightarrow \mathbb{R}^r/D$ be the (equivariant) map obtained by composing the map $F$ with the projection $\mathbb{R}^r \rightarrow \mathbb{R}^r/D_{\mathbb{R}^r}$, where $D_{V^r}$ stands for the diagonal $D_{V^r} \subset V^r$. This map must have a zero.

 \medskip
 In the opposite case there arises an equivariant map $\mathcal{C}_2 \rightarrow S^{r-2}$.
This is,  in light of the $(r-2)$-connectivity of $\mathcal{C}_2$  \cite{ZV92}, a contradiction with either Dold's theorem \cite{Mat} or (more general) Volovikov's Theorem \cite{Vol96-1}.

It remains to apply D. Gale's original trick \cite{G}, already described in the proof of Theorem \ref{ThmPrime}.
 \qed

\medskip\noindent
\emph{New proof of Theorem \ref{ThmKA}.}
Based on  preferences $(A_i^j)$, let us construct  new   preferences $(C_i^j)=(\phi (A_i^j))$ defined on the configuration space $\mathcal{C}_2$.

\medskip

The map $\phi$ can be described in  a form of an algorithm:

\begin{enumerate}
\item We are given  a cut $x$ of the segment $I$ with $2r-2$ cut points.  Create an induced cut $y$ of $I$ with $r-1$ cut points  ($y\in \Delta^{r-1}$), preserving the same collection of non-degenerate intervals $\{I^\nu_k\}_{k=1}^s$. In other words we eliminate $r-1$ superfluous (multiple) cuts. Note that this step can be performed in many different ways.

\item If the preferences $\{A_i^j\}$ of a player dictate the choice of some non-degenerate tile $i$,   add the box $\alpha(i)$ to the preferences of the player.
 Note that the property $\mathbf{P_{pe}}$ of $(A_i^j)$) implies that no matter how the superfluous cuts are eliminated, the result will be one and the same.
 \item If the preferences $(A_i^j)$ dictate a player to choose a degenerate interval (which occurred  after $r-1$ cuts), observe
 that there necessarily exists an empty box, since in this case the number of non-degenerate tiles is at most $r-1$. Add all empty boxes to the preferences of the player.
\end{enumerate}

\medskip
We emphasize once again that the condition $\mathbf{P_{pe}}$ is used in the proof that the preferences $C_i^j$ are well-defined.

\medskip
It is not difficult to check that the preferences $(C_i^j)=(\phi(A_i^j))$ satisfy all conditions of Theorem \ref{ThmPrimePower}. So there exist an envy-free division with $(C_i^j)$, and, by construction, a full envy-free division for $(A_i^j)$.
 \qed

\section{Configuration space $[r]^{*r}$}\label{sec:ManyPieces}

We illustrate the versatility of the method  of configuration spaces by one more example.
 As in Section \ref{Sec1} we consider partitions of the interval $I=[0,1]$ into $r$ tiles by $r-1$ cuts.  However in this section we put no restriction on the allocation functions $\alpha : [r]\rightarrow [r]$ and consider all possible ways to allocate $r$ tiles to $r$ boxes. In particular, a box may be empty or contain several tiles.

\medskip

The corresponding configuration space  $\mathcal{C}_3 = [r]^{*r}$ is the $r$-fold join of the $0$-dimensional simplicial complex $[r]$. Recall \cite{Mat, Z17} that this configuration space naturally arises in the topological Tverberg problem (a monochromatic predecessor of the colored Tverberg problem).

\begin{thm}\label{ThmPrimePower2}
Let $r$ be a prime power.
  Let $$(D_i^j)_{i,j=1}^r,\ \  D_i^j\subseteq \mathcal{C}_3$$ be a matrix of the new-style (closed, covering, equivariant) preferences. Then there exists an envy-free division:
\begin{equation}\label{eqn:empty-5}
     \bigcap_{i\in [r]}  D_i^i \neq\emptyset \, .
\end{equation}
In other words there exists a cut and an allocation of non-degenerate tiles into the boxes, such that the player $j$ prefers the box $j$ for all $ j$.
\end{thm}
\proof

The proof follows the pattern of the proof of Theorem \ref{ThmPrime}.  Assuming the contrary, we again arrive at an equivariant map $\mathcal{C}_3\rightarrow S^{r-2}$ from the configuration space to the sphere.
The non-existence of such a map is a well-known fact \cite{Mat, Z17}, as one of the key ingredients  in the proof of the topological  Tverberg theorem.
 \qed

\section{Appendix. Chessboard complexes}\label{sec:chessboard}
\label{sec:appendix}

The \emph{Chessboard complex}   $\Delta_{m,n}$ is  the simplicial complex of all non-attacking placements  of rooks on a $[m]\times [n]$ chessboard. More precisely the vertex-set of $\Delta_{m,n}$ is the set $[m]\times [n]$ (interpreted as a $(m\times n)$-chessboard) and $S\subset [m]\times [n]$ is a simplex if $S$ has at most one element in each row or column of $[m]\times [n]$.

\medskip
Similarly the join $\Delta_{r,r-1}*[r]$, introduced in Section \ref{Sec1}, is the simplicial complex whose vertices $v_{ij}$ correspond to the cells of the chessboard $r\times r$, and a set of vertices
corresponds to a simplex whenever (1) there is at most one vertex in each column, and (2) there is at most one vertex in each row, with a unique possible exception, when there are two vertices in a row, and one of these two vertices belongs to the column $r$.
From this description it is clear that the rows of the chessboard correspond to the boxes, whereas the columns correspond to the tiles of a cut.

\medskip
If a column of the chessboard $[r]\times [r]$ is empty, the corresponding tile degenerates.
Each point $\hat{z}\in \Delta_{r,r-1}*[r]$ is in (the interior of) some simplex, and has a (unique) representation $\hat{z} = \sum_{i=1}^r z_i v_{i,\alpha(i)}$, in terms of the corresponding barycentric coordinates.
This point is associated to
(1) the cut whose $z$-representation (Section \ref{sec:conf-old-new}) is $(z_1,...,z_r)$, and (2) the allocation function $\alpha$ which allocates the tile $i$ to the corresponding box $\alpha(i)$.

In other words, the join  $\Delta_{r,r-1}*[r]$ is precisely the configuration space which parameterizes all partitions/allocations described in Section \ref{Sec1}.
\medskip

Analogously, the complex $\Delta_{r,2r-1}$ is isomorphic to the configuration space from Section \ref{Sec2} which records partitions into $2r-1$ tiles with at most $r$ non-degenerate ones and the corresponding allocation functions.

\medskip

The $(r-2)$-connectedness of the complex $\Delta_{r,2r-1}$ was used in \cite{ZV92} for the proof of the original {\em colored Tverberg theorem}. The complex $\Delta_{r,r-1}$ is always a pseudomanifold and its properties were used in the proof of the \emph{optimal colored Tverberg theorem} \cite{BMZ} (see also \cite{vz11, jpz-3}). The role of chessboard complexes and other configuration spaces used in combinatorics and discrete geometry is discussed in \cite{Mat, Z17}.

\subsection*{Acknowledgements} We are grateful to the anonymous referee of the paper \cite{jpz-2}  who pointed out that our methods imply Theorem \ref{ThmKA}.
 R. \v Zivaljevi\' c was supported by the Serbian Ministry of Education, Science and Technological Development through Mathematical Institute of the Serbian Academy of Sciences and Arts. Sections 3 and 5  are  supported by the Russian Science Foundation under grant  21-11-00040.

\end{document}